\documentclass[12pt]{article}
\usepackage[utf8]{inputenc}
\usepackage[english]{babel}
\usepackage{graphicx}
\usepackage{vmumath}
\usepackage{amsmath}
\usepackage{amsthm}
\usepackage{avo_book}
\usepackage{tikz}
\usepackage{caption}
\usepackage{subcaption}
\usepackage{booktabs}
\usepackage[a4paper,left=15mm,right=15mm,top=30mm,bottom=20mm]{geometry}
\noaftermath
\begin{document}

\bibliographystyle{unsrt}

\title{Unsensed enumeration of cubic unicellular maps on orientable and non‑orientable surfaces: explicit formulas via orbifolds}

\author{
Alexander Omelchenko\thanks{Constructor University Bremen, Campus Ring 1, 28759 Bremen, Germany. E-mail: {\tt aomelchenko@constructor.university}} \and
Igor Labutin\thanks{Constructor University Bremen, Campus Ring 1, 28759 Bremen, Germany. E-mail: {\tt igor.labutin@gmail.com}}
}

\begin{abstract}
We enumerate cubic (3‑regular) unicellular maps on closed surfaces up to all homeomorphisms. Using the orbifold approach, we reduce the unsensed enumeration to explicit counts of quotient maps and rooted cubic/precubic maps on simpler surfaces. For orientable hosts this yields a compact identity expressed through known sensed and rooted numbers; for non‑orientable hosts we obtain a fully explicit finite‑sum expression via precubic counts. Numerical tables are provided, together with a brief asymptotic discussion.

\end{abstract}

\maketitle

\section{Introduction}

A \emph{one-face} (or \emph{unicellular}) topological map $M$ on a closed surface $X$ is a 2-cell embedding of a connected graph $G$ (loops and multiple edges allowed) into a compact connected two-dimensional manifold without boundary such that the complement $X - G$ is a single open 2-cell; the vertices and edges of $M$ are its $0$- and $1$-cells, respectively \cite{Handbook_of_Graph_Theory}.  
We consider both orientable and non-orientable surfaces. Each such surface is determined by its genus $g$: an orientable surface of genus $g$ is a sphere with $g$ handles, while a non-orientable surface of genus $g$ is a sphere with $g$ crosscaps. The Euler characteristic is $\chi = 2 - 2g$ for orientable surfaces $X^+$ and $\chi = 2 - g$ for non-orientable ones $X^-$.  

Two maps on the same surface are \emph{isomorphic} if a homeomorphism of the surface induces an isomorphism of their underlying graphs. The set of all maps thus splits into unlabelled isomorphism classes. On orientable surfaces one distinguishes between orientation-preserving and orientation-reversing homeomorphisms: classes modulo the former are called \emph{sensed maps}, while classes modulo all homeomorphisms (including reversing ones) are the \emph{unsensed maps}. Enumeration of unsensed maps is a substantially harder problem, as it requires accounting for orientation-reversing symmetries.  

The enumeration of sensed maps was initiated by Liskovets for planar maps \cite{Liskovets_85} and extended to arbitrary orientable surfaces by Mednykh and Nedela \cite{Mednykh_Nedela}. Their approach, now classical, reduces counting to enumerating \emph{quotient maps on orbifolds} under finite cyclic group actions. This orbifold framework has been further developed for hypermaps \cite{Mednykh_Hypermaps}, for one-face families of higher valence \cite{Krasko_Omelch_4_reg_one_face_maps}, and for regular maps on the torus and higher-genus surfaces \cite{Torus_Part_I,POMI_Reg_maps_english}.  
On the unsensed side, orientable maps of all genera were treated in \cite{Azevedo}, explicit formulas for unsensed $r$-regular maps on the torus were obtained in \cite{Torus_Part_II}, and a general orbifold scheme for unsensed maps on both $X_g^+$ and $X_g^-$ was formulated in \cite{Unsensed_Maps}.  

Recent developments continue to explore enumeration of maps and hypermaps with prescribed symmetries and automorphism groups on higher surfaces, see e.g. \cite{Chapuy_Do_Fang_2021,Kang_Nedela_2021,Mednykh_Nedela_Stukachev_2023}.  
In contrast to these general frameworks, the present paper provides fully explicit closed formulas for a specific regular class — cubic (3-regular) unicellular maps — thereby completing the enumeration program for unsensed regular maps initiated in \cite{Unsensed_Maps}.  

For orientable hosts, every orientation-reversing symmetry yields an orbifold with a single boundary component; for non-orientable hosts we obtain either a boundary orbifold with controlled index-2 branch points or a closed non-orientable orbifold with branch indices among $\{2,3,\ell\}$ and exactly one index-$\ell$ point in the unique face.  
These constraints reduce the problem to explicit counts of rooted cubic and precubic maps on simpler closed surfaces. The required building blocks are classical: Walsh–Lehman formulas for orientable unicellular maps \cite{Walsh_Lehman} and Bernardi–Chapuy expressions for non-orientable precubic maps \cite{Chapuy_non_orientable}; planar precursors go back to Tutte’s seminal enumerations \cite{Tutte_Census,Tutte_triangulations}.  

By combining the orbifold reduction with these rooted and precubic counts, we derive fully explicit closed formulas for the numbers of unsensed cubic one-face maps on both orientable and non-orientable surfaces. Numerical tables up to genus $20$ and a short asymptotic discussion are also provided.  
This paper completes the unsensed enumeration program for regular unicellular maps, providing the first fully closed expressions for the cubic case.

\section{Enumeration of unsensed $3$-regular one-face maps on orientable surfaces}

Let $X_g^+$ be a closed orientable surface of genus $g$. Following the general orbifold approach of Mednykh and Nedela \cite{Mednykh_Nedela}, the number $\tilde{\tau}_{X_g^+}(n)$ of \emph{sensed} maps with $n$ edges is given by
\begin{equation}
\label{eq:Mednykh_Ned_2006_new}
\tilde{\tau}_{X_g^+}(n)=\frac{1}{2n}\sum_{\substack{\ell\mid 2n\\ \ell\,m=2n}}\ \sum_{O\in {\rm Orb}(X_g^+/\mathbb{Z}_{\ell})}
\Epi_o(\pi_1(O), \mathbb{Z}_{\ell})\ \tau_O(m).
\end{equation}
Here $O=X_g^+/\mathbb{Z}_\ell$ is a quotient orbifold of $X_g^+$, $\Epi_o(\pi_1(O), \mathbb{Z}_{\ell})$ is the number of order‑preserving epimorphisms, and $\tau_O(m)$ is the number of rooted quotient maps on $O$ with $m$ darts. For our purposes we only need that \eqref{eq:Mednykh_Ned_2006_new} provides the sensed counts as an \emph{input}; in particular, for $3$‑regular one‑face maps on $X_g^+$ the closed form for $\tilde\tau_+^{(3)}(g)$ is available in \cite[formula (20)]{POMI_Reg_maps_english}.

For \emph{unsensed} maps on $X_g^+$ the corresponding formula (see \cite{Unsensed_Maps}) has the form
\begin{equation}
\label{eq:unsensed_orientable_final_new}
\bar\tau_{X_g^+}(n) = \frac{1}{2}\left( \tilde\tau_{X_g^+}(n) + \frac{1}{2n} \sum_{\substack{m\mid 2n\\ \ell m = n}} \ \sum_{O \in {\rm Orb}^-(X_g^+/\mathbb{Z}_{2\ell})} \tau_O(2m)\, \Epi_o^+(\pi_1(O), \mathbb{Z}_{2\ell}) \right),
\end{equation}
where ${\rm Orb}^-(X_g^+/\mathbb{Z}_{2\ell})$ denotes the orbifolds arising from orientation‑reversing homeomorphisms, and $\Epi_o^+$ counts orientation‑ and order‑preserving epimorphisms.

For cubic one‑face maps on $X_g^+$ we have, by Euler and the handshaking lemma,
\[
n=6g-3,\qquad k=4g-2.
\]
The contribution in \eqref{eq:unsensed_orientable_final_new} coming from orientation‑reversing symmetries simplifies drastically; the key facts are the following.

\begin{propos_eng}\label{prop:boundary_no_branches}
For a $3$‑regular one‑face map on $X_g^+$, any orbifold $O$ corresponding to an orientation‑reversing homeomorphism is a surface (orientable or non‑orientable) with boundary and with no branch points.
\end{propos_eng}

\evidpEng The statement can be seen in two ways. First, one‑face maps on orientable surfaces admit a representation by chord diagrams; any symmetry of such a diagram is either a rotation or a reflection. Reflections correspond to orientation‑reversing homeomorphisms of period $2$, hence boundary must appear in the quotient, and no branch point can lie in the unique face. Second, in the general orbifold language, an orientation‑reversing symmetry of $X_g^+$ produces either a non‑orientable closed orbifold or an orbifold with boundary. In the closed non‑orientable case, to keep one face after lifting, a branch point of index $2\ell$ must lie in the face; however, for the coefficients in \eqref{eq:unsensed_orientable_final_new} one has $\Epi_o^+(\pi_1(O),\mathbb{Z}_{2\ell})=0$ for such $O$ (see \cite{Unsensed_Maps}). In the boundary case, period is necessarily $2$, and all branch indices (if any) equal $2$; but orbifolds with boundary cannot carry index‑$2$ branch points simultaneously (cf.\ \cite[§2]{Azevedo}). Hence only boundary remains and no branch points occur. \qed

\begin{conseq_eng}\label{cor:l2}
In \eqref{eq:unsensed_orientable_final_new} only the case $\ell=1$ contributes, so
\begin{equation}
\label{eq:unsensed_orientable_final_l2_new}
\bar\tau_{X_g^+}(n) = \frac{1}{2}\left( \tilde\tau_{X_g^+}(n) + \frac{1}{2n}\sum_{O\in {\rm Orb}^-(X_g^+/\mathbb{Z}_{2})} \tau_O(2n)\right).
\end{equation}
\end{conseq_eng}

In the cubic case the boundary is in fact unique.

\begin{figure}[t]
\centering
\includegraphics[scale=0.75]{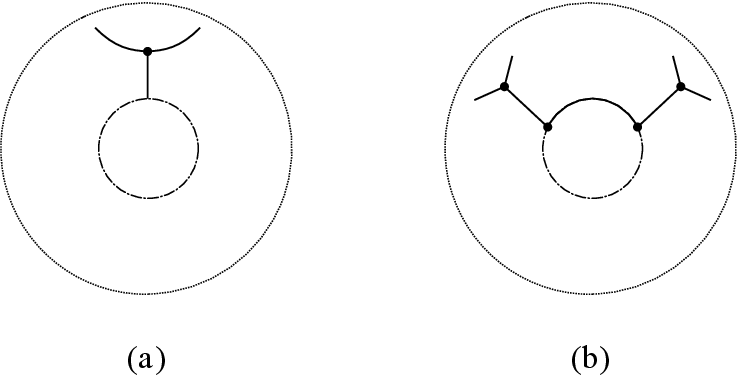}
\caption{Boundary configurations for quotient maps on $O$: (a) a half‑edge ending on the boundary; (b) a boundary edge with two incident cubic vertices, each incident with exactly one normal edge.}
\label{fig:Boundary_example}
\end{figure}

\begin{propos_eng}\label{prop:single_boundary}
For $3$‑regular one‑face maps on $X_g^+$, the orbifold $O$ corresponding to an orientation‑reversing homeomorphism has a \emph{single} boundary component.
\end{propos_eng}

\evidpEng Along the boundary only the two local configurations in Fig.~\ref{fig:Boundary_example} may occur. In both, a boundary component cannot be covered entirely by quotient edges, hence each boundary component meets the unique face. If there were two boundary components (or the same component met the face more than once), lifting would create a face non‑homeomorphic to a disk, contradicting unicellularity. \qed

Now apply the Riemann–Hurwitz relation
\begin{equation}
\label{eq:Riemann_Hurwitz_new}
-\chi=\ell\left(\alpha\,\mathfrak g-2+h+\sum_{i=1}^r\Bigl(1-\frac{1}{m_i}\Bigr)\right),
\end{equation}
with $\ell=2$, $h=1$, $r=0$. This gives the precise form of $O$:

\begin{propos_eng}\label{prop:genus_orbifold}
If $g$ is odd, then $O$ is non‑orientable of genus $\mathfrak g=g$. If $g$ is even, then either $O$ is orientable of genus $\mathfrak g=g/2$, or $O$ is non‑orientable of genus $\mathfrak g=g$.
\end{propos_eng}

The enumeration of quotient maps on $O$ now reduces to rooted cubic maps on closed surfaces.

\begin{propos_eng}\label{prop:reduction_rooted}
Contract the unique boundary component of a rooted quotient map on $O$. One obtains a rooted map on a closed surface with one degree‑$2$ vertex; contracting this vertex yields a rooted $3$‑regular one‑face map with
\[
n'=\frac{n-3}{2}
\]
edges on $X_{g/2}^+$ (when $g$ is even and $O$ is orientable) or on $X_g^-$ (when $O$ is non‑orientable). Conversely, each rooted cubic one‑face map of the corresponding type yields exactly $2n$ rooted quotient maps on $O$.
\end{propos_eng}

Therefore
\[
\tau_{O^+}(2n)=2n\cdot \tau_+^{(3)}(g/2)\quad (g\ \text{even}),\qquad
\tau_{O^-}(2n)=2n\cdot \tau_-^{(3)}(g),
\]
where
\begin{equation}
\label{eq:num_rooted_3_orient_new}
\tau_+^{(3)}(g)=\frac{2\,(6g-3)!}{12^g\,g!\,(3g-2)!}
\end{equation}
and
\begin{equation}
\label{eq:num_rooted_3_non_orient_new}
\tau_-^{(3)}(g)=
\begin{cases}
\displaystyle
\frac{2^{\,2h-2}\,h!\,(6h-2)!}{3^{\,h-1}\,(2h)!\,(3h-1)!}\,\sum_{i=0}^{h-1}\binom{2i}{i}\,16^{-i}, & g=2h,\\[1.2ex]
\displaystyle \frac{2^{\,6h}\,(3h)!}{3^{\,h}\,h!}, & g=2h+1,
\end{cases}
\end{equation}
see \cite{Walsh_Lehman,Chapuy_non_orientable}.

Putting these expressions into \eqref{eq:unsensed_orientable_final_l2_new} we obtain the final identity.

\begin{theor_eng}\label{thm:orientable_final}
For $g\ge1$,
\begin{equation}
\label{eq:num_3_one_faced_unsensed_or_map_new}
\bar\tau_+^{(3)}(g)=\frac{1}{2}\left(\tilde\tau_+^{(3)}(g)\ +\ \tau_+^{(3)}(g/2)\ +\ \tau_-^{(3)}(g)\right),
\end{equation}
with the convention $\tau_+^{(3)}(g/2)=0$ for odd $g$. Here $\tilde\tau_+^{(3)}(g)$ is taken from \cite{POMI_Reg_maps_english}.
\end{theor_eng}

\section{Enumeration of unsensed $3$-regular one-face maps on non-orientable surfaces}

Let $X_g^-$ be a closed non-orientable surface of genus $g$. For a cubic one-face map on $X_g^-$ we have $n=3g-3$ edges and $k=2g-2$ vertices by Euler’s formula and the handshaking lemma. The unsensed numbers are given by the general orbifold formula (see \cite{Unsensed_Maps})
\begin{equation}
\label{eq:unsensed_nonorientable_master}
\bar\tau_{X_g^-}(n)=\frac{1}{4n}\sum_{\substack{m\mid 2n\\ \ell m=2n}}\ \sum_{O\in {\rm Orb}(X_g^-/\mathbb{Z}_\ell)} \tau_O(2m)\,\bigl(\Epi_o(\pi_1(O),\mathbb{Z}_\ell)-\Epi_o^+(\pi_1(O),\mathbb{Z}_\ell)\bigr),
\end{equation}
where $\tau_O$ counts rooted quotient maps on the orbifold $O$ with the prescribed number of darts, and the multiplicity coefficients are the differences of (order‑preserving) epimorphism numbers. This is our starting point for the cubic unicellular case.  

\medskip
We first describe the relevant orbifolds. There are two families. In the \emph{boundary} family the period is $\ell=2$ and the quotient $O$ has one boundary component; index‑$2$ branch points (if any) may lie only at vertices of the quotient map or at free ends of semi‑edges; no branch point can lie in the unique face. In the \emph{closed} family the quotient $O$ is a closed non‑orientable orbifold with branch indices among $\{2,3,\ell\}$ and with exactly one index‑$\ell$ point inside the unique face; the remaining branch points coincide with vertices (index~3) or with free ends of semi‑edges (index~2). Both statements are standard in this setting and will be used without further reference.

\begin{propos_eng}\label{prop:RH-nonori}
In the boundary case $(\ell=2)$, by the Riemann–Hurwitz relation with one boundary component one has
\[
g=\begin{cases}
4\mathfrak g+r,& O\ \text{orientable},\\
2\mathfrak g+r,& O\ \text{non-orientable},
\end{cases}
\]
so that $r=g-4\mathfrak g$ in the orientable case (hence $0\le \mathfrak g\le\lfloor g/4\rfloor$) and $r=g-2\mathfrak g$ in the non‑orientable case (hence $1\le \mathfrak g\le\lfloor g/2\rfloor$). In the closed case $(\ell>2)$ the orbifold signatures have the form
\[
O^-(\mathfrak g;\ 2^{n_s},\,3^{n_v},\,\ell),\qquad 6g-6=\ell\,(6\mathfrak g-6+3n_s+4n_v),
\]
with the divisibility constraints $2\mid \ell\Rightarrow n_s>0$ and $3\mid \ell\Rightarrow n_v>0$.
\end{propos_eng}

\evidpEng This is a direct application of the Riemann–Hurwitz formula to orbifolds with one boundary component (for $\ell=2$) and to closed non‑orientable orbifolds (for $\ell>2$), keeping track of local lifting rules for semi‑edges, boundary edges and vertex lifts.  

\medskip
The multiplicity coefficients in \eqref{eq:unsensed_nonorientable_master} are explicit. In the boundary family ($\ell=2$) they reduce to simple powers of $2$.

\begin{propos_eng}[Epimorphisms for $\ell=2$]\label{prop:epi-l2}
Let $O^+$ (resp.\ $O^-$) be a boundary orbifold of genus $\mathfrak g$ that is orientable (resp.\ non‑orientable), with $r$ index‑2 branch points. Then
\[
\bigl(\Epi_o-\Epi_o^+\bigr)(\pi_1(O^+),\mathbb{Z}_2)=
\begin{cases}
2^{2\mathfrak g},& r>0,\\[0.2ex]
2^{2\mathfrak g}-1,& r=0,
\end{cases}
\qquad
\bigl(\Epi_o-\Epi_o^+\bigr)(\pi_1(O^-),\mathbb{Z}_2)=
\begin{cases}
2^{\mathfrak g},& r>0,\\[0.2ex]
2^{\mathfrak g}-1,& r=0.
\end{cases}
\]
\end{propos_eng}

\evidpEng This follows from the explicit formulas for epimorphisms from orbifold fundamental groups with boundary (Jordan totients), with the orientation‑preserving term present only when $r=0$.  

\medskip
The quotient‑to‑rooted reduction in the boundary family mirrors the orientable case. Contracting the unique boundary component produces a rooted map on a closed surface with one degree‑$2$ vertex; removing this vertex yields a rooted precubic map with a prescribed number of leaves.

\begin{propos_eng}[Reduction to precubic counts]\label{prop:precubic-red}
In the boundary case, quotient maps on $O$ with $2n$ darts reduce bijectively (with the usual root normalisation) to rooted precubic one‑face maps on closed surfaces of genus $\mathfrak g$, with exactly $r$ leaves, where $(\mathfrak g,r)$ are as in Proposition~\ref{prop:RH-nonori}. Quantitatively,
\[
\tau_O(2n)=2n\cdot 
\begin{cases}
\tau^{(1\!\div\!3)}_+(\mathfrak g,r),& O=O^+ \text{ (orientable)},\\
\tau^{(1\!\div\!3)}_-(\mathfrak g,r),& O=O^- \text{ (non‑orientable)}.
\end{cases}
\]
\end{propos_eng}

\evidpEng The $2n$ factor is the number of admissible root positions that survive the contraction/expansion moves; leaves cannot serve as roots.  

\medskip
The needed precubic inputs are classical. For orientable hosts (Walsh–Lehman),
\begin{equation}
\label{eq:precubic-or}
\tau^{(1\!\div\!3)}_+(\mathfrak g,r)=\frac{2\,(2r+6\mathfrak g-3)!}{12^{\mathfrak g}\,\mathfrak g!\,(r+3\mathfrak g-2)!\,r!}.
\end{equation}
For non‑orientable hosts (Bernardi–Chapuy),
\begin{equation}
\label{eq:precubic-nor}
\tau^{(1\!\div\!3)}_-(\mathfrak g,r)=
\begin{cases}
\displaystyle \frac{2\,c_h\,(2r+6h-3)!}{r!\,(r+3h-2)!}, & \mathfrak g=2h,\\[1.2ex]
\displaystyle \frac{2^{\,6h+2r}}{3^h\,h!\,(r+3h)!}, & \mathfrak g=2h+1,
\end{cases}
\qquad
c_h=\frac{2^{\,2h-2}h!}{3^{\,h-1}(2h)!}\sum_{i=0}^{h-1}\binom{2i}{i}16^{-i}.
\end{equation}
Substituting \eqref{eq:precubic-or}–\eqref{eq:precubic-nor} and Proposition~\ref{prop:epi-l2} into \eqref{eq:unsensed_nonorientable_master} with $\ell=2$ gives the boundary contribution
\begin{equation}
\label{eq:l2-contrib}
\bar\tau^{(3)}_-(g)\big|_{\ell=2}=\frac12\sum_{\mathfrak g=0}^{\lfloor g/4\rfloor}\!\!\bigl(2^{2\mathfrak g}-\mathbf{1}_{\{r=0\}}\bigr)\,\tau^{(1\!\div\!3)}_+(\mathfrak g,\,g-4\mathfrak g)\ +\ \frac12\sum_{\mathfrak g=1}^{\lfloor g/2\rfloor}\!\!\bigl(2^{\mathfrak g}-\mathbf{1}_{\{r=0\}}\bigr)\,\tau^{(1\!\div\!3)}_-(\mathfrak g,\,g-2\mathfrak g),
\end{equation}
where in each sum $r$ is the corresponding number of index‑$2$ branch points (so $r=0$ occurs only when $g=4\mathfrak g$ in the first sum, and when $g=2\mathfrak g$ in the second).  

\medskip
We now turn to the closed family $(\ell>2)$. The admissible orbifolds are precisely those with signatures $O^-(\mathfrak g;2^{n_s},3^{n_v},\ell)$ obeying $6g-6=\ell(6\mathfrak g-6+3n_s+4n_v)$ and the divisibility constraints. For these,
\begin{equation}
\label{eq:epi-l>2}
\bigl(\Epi_o-\Epi_o^+\bigr)(\pi_1(O^-),\mathbb{Z}_\ell)=\varepsilon(\ell,\mathfrak g,n_v)=
\begin{cases}
\ell^{\mathfrak g-1}\,\varphi(\ell)\,2^{n_v},& \ell\ \text{odd},\\[0.4ex]
2\,\ell^{\mathfrak g-1}\,\varphi(\ell)\,2^{n_v},& \ell\ \text{even and }\ \frac{\ell}{2}n_s+\frac{\ell}{3}n_v+1\ \text{even},\\[0.4ex]
0,& \text{otherwise}.
\end{cases}
\end{equation}
The quotient‑map enumeration reduces to precubic maps with $k=n_s+n_v$ leaves on $X_{\mathfrak g}^-$; one must additionally choose which $n_s$ of the $k$ leaves arise from semi‑edges (binomial factor $\binom{k}{n_s}$) and adjust the root normalisation by the ratio of dart counts:
\[
\frac{\#\text{darts of quotient}}{\#\text{darts of precubic}}=\frac{(6\mathfrak g-6)/\ell}{(6\mathfrak g-6)/\ell+n_s}.
\]
Putting this into \eqref{eq:unsensed_nonorientable_master} yields a finite sum over all admissible tuples $(\ell,\mathfrak g,n_s,n_v)$:
\begin{equation}
\label{eq:l>2-contrib}
\bar\tau^{(3)}_-(g)\big|_{\ell>2}=\frac{1}{4}\sum_{(\ell,\mathfrak g,n_s,n_v)}\varepsilon(\ell,\mathfrak g,n_v)\,\binom{n_s+n_v}{n_s}\,\frac{(6\mathfrak g-6)/\ell}{(6\mathfrak g-6)/\ell+n_s}\ \tau^{(1\!\div\!3)}_-(\mathfrak g,n_s+n_v),
\end{equation}
equivalently, after cancelling the outer factor $1/(4n)$ with $n=3g-3$, the factor $\frac{(6\mathfrak g-6)/\ell}{(6\mathfrak g-6)/\ell+n_s}$ may be written as $\bigl(3g-3+\tfrac{\ell}{2}\,n_s\bigr)^{-1}$.  

\medskip
Finally, the identity term in \eqref{eq:unsensed_nonorientable_master} contributes $\frac{1}{4n}\tau_-^{(3)}(g)=\frac{1}{4(3g-3)}\tau_-^{(3)}(g)$. Summarising the three parts we obtain the explicit formula for the non‑orientable host.

\begin{theor_eng}\label{thm:nonori-final}
For $g\ge 2$,
\begin{equation}
\label{eq:final-nonori}
\bar\tau^{(3)}_-(g)=\frac{1}{4(3g-3)}\,\tau_-^{(3)}(g)\ +\ \bar\tau^{(3)}_-(g)\big|_{\ell=2}\ +\ \bar\tau^{(3)}_-(g)\big|_{\ell>2},
\end{equation}
with $\bar\tau^{(3)}_-(g)\big|_{\ell=2}$ as in \eqref{eq:l2-contrib} and $\bar\tau^{(3)}_-(g)\big|_{\ell>2}$ as in \eqref{eq:l>2-contrib}. The rooted cubic inputs are
\begin{equation}
\label{eq:rooted-cubic-nonori}
\tau_-^{(3)}(g)=
\begin{cases}
\displaystyle \frac{2^{\,2h-2}\,h!\,(6h-2)!}{3^{\,h-1}\,(2h)!\,(3h-1)!}\,\sum_{i=0}^{h-1}\binom{2i}{i}16^{-i}, & g=2h,\\[1.2ex]
\displaystyle \frac{2^{\,6h}\,(3h)!}{3^{\,h}\,h!}, & g=2h+1,
\end{cases}
\end{equation}
and the precubic inputs are given by \eqref{eq:precubic-or}–\eqref{eq:precubic-nor}.
\end{theor_eng}

\noindent
All coefficients and summation ranges are explicit, so \eqref{eq:final-nonori} is a completely explicit closed formula; its first values agree with the numerical tables included later in the paper.

\section{Numerical data and asymptotics}

We collect numerical values for cubic one-face maps on $X_g^{\pm}$ and record simple asymptotic consequences of the explicit formulas proved above. For orientable hosts $X_g^+$ we list rooted, sensed and unsensed numbers; for non‑orientable hosts $X_g^-$ we list rooted and unsensed numbers. All values agree with the closed forms stated in Theorems~\ref{thm:orientable_final} and~\ref{thm:nonori-final}.

{\setlength{\tabcolsep}{3pt}
\begin{table}[t!]
\begin{center}
\footnotesize
\begin{tabular}{cccc}
\midrule
$g$ & $\tau_{+}^{(3)}(g)$ & $\tilde\tau_{+}^{(3)}(g)$ & $\bar\tau_{+}^{(3)}(g)$ \\ 
\midrule
1 & 1 & 1 & 1 \\
2 & 105 & 9 & 8 \\
3 & 50050 & 1726 & 927 \\
4 & 56581525 & 1349005 & 676445 \\
5 & 117123756750 & 2169056374 & 1084610107 \\
6 & 386078943500250 & 5849686966988 & 2924847922929 \\
7 & 1857039718236202500 & 23808202021448662 & 11904101304325611 \\
8 & 12277353837189093778125 & 136415042681045401661 & 68207521363461659373 \\
9 & 106815706684397824557193750 & 1047212810636411989605202 & 523606405320272947813801 \\
10 & 1183197582943074702620035168750 & 10378926166167927379808819918 & 5189463083084174721816125584 \\
\midrule
\end{tabular}
\caption{Cubic one-face maps on orientable surfaces $X_g^{+}$}
\label{table:orientable_surface}
\end{center}
\end{table}}

{\setlength{\tabcolsep}{3pt}
\begin{table}[t!]
\begin{center}
\footnotesize
\begin{tabular}{ccc}
\midrule
$g$ & $\tau_{-}^{(3)}(g)$ & $\bar\tau_{-}^{(3)}(g)$ \\ 
\midrule
2 & 6 & 2 \\
3 & 128 & 11 \\
4 & 3780 & 144 \\
5 & 163840 & 3627 \\
6 & 8828820 & 149288 \\
7 & 587202560 & 8170800 \\
8 & 45821335560 & 545671762 \\
9 & 4133906022400 & 43063046307 \\
10 & 421946699674500 & 3906934079662 \\
11 & 48151737348915200 & 401264673924438 \\
12 & 6070544859205827000 & 45988979036528440 \\
13 & 838225443769915801600 & 5821010056777072838 \\
14 & 125787689149526729325000 & 806331341176441101980 \\
15 & 20385642792484352294912000 & 121343111865634574938768 \\
16 & 3548258423062128985899690000 & 19712546794881999409462482 \\
17 & 660168656191813264718430208000 & 3438378417666873290074260643 \\
18 & 130746565669943973430227429382500 & 640914537597785062325259175158 \\
19 & 27463016097579431812286696652800000 & 127143593044349500804170430994988 \\
20 & 6098023559259606741021710317037175000 & 26745717365173718867249062116990380 \\
\midrule
\end{tabular}
\caption{Cubic one-face maps on non‑orientable surfaces $X_g^{-}$}
\label{table:non_orientable_surface}
\end{center}
\end{table}}

For orientable hosts, Stirling’s approximation applied to the rooted formula
\[
\tau^{(3)}_+(g)=\frac{2\,(6g-3)!}{12^g\,g!\,(3g-2)!}
\]
gives the growth
\[
\tau^{(3)}_+(g)=\Theta\!\left(g^{\,2g-\frac{3}{2}}\left(\frac{144}{e^2}\right)^{\!g}\right),
\]
and in the sensed master sum the identity term dominates: the total contribution of $L\in\{2,3,6\}$ is exponentially small in $g$. Consequently,
\[
\tilde\tau^{(3)}_+(g)=\frac{\tau^{(3)}_+(g)}{2E}\,(1+o(1)),\qquad
\bar\tau^{(3)}_+(g)=\frac{\tau^{(3)}_+(g)}{4E}\,(1+o(1)),\qquad E=6g-3.
\]
In the non‑orientable case, from the master formula and Theorem~\ref{thm:nonori-final} the identity contribution equals $\tau_-^{(3)}(g)/(4n)$ with $n=3g-3$, while the remaining terms have strictly smaller exponential order; numerically one observes
\[
\bar\tau^{(3)}_-(g)=\frac{\tau^{(3)}_-(g)}{4n}\,(1+o(1)).
\]
We also checked the first values by independent map generation; the results coincide with the tables above.

\section*{Conclusion}

We obtained explicit formulas for the numbers of unsensed cubic one‑face maps on orientable and on non‑orientable surfaces of genus $g$. On $X_g^+$ the boundary structure of the orbifold under orientation‑reversing symmetries reduces the fixed‑map enumeration to rooted counts, which yields a compact identity linking unsensed, sensed and rooted numbers. On $X_g^-$ the admissible orbifold signatures are completely described, the epimorphism coefficients are explicit, and quotient maps reduce to precubic counts, leading to a finite closed sum. Numerical tables and simple asymptotic rules are presented. 

The authors thank E. Krasko for helpful discussions and constructive remarks.


\newpage

\begin{thebibliography}{10}

\bibitem{Azevedo}
A. Breda d'Azevedo, A. Mednykh and R. Nedela.
\newblock Enumeration of maps regardless of genus: Geometric approach.
\newblock {\em Discrete Mathematics}, 310:1184--1203, 2010.

\bibitem{Chapuy_non_orientable}
O. Bernardi, G. Chapuy.
\newblock Counting unicellular maps on non-orientable surfaces.
\newblock {\em Advances in Applied Mathematics}, 48(2):259--275, 2011.

\bibitem{Chapuy}
G. Chapuy.
\newblock A new combinatorial identity for unicellular maps, via a direct bijective approach.
\newblock {\em Advances in Applied Mathematics}, 47(4):874--893, 2011.

\bibitem{Chapuy_Do_Fang_2021}
G. Chapuy, M. Do, L. Fang.
\newblock Unicellular maps, topological recursion and Hurwitz numbers.
\newblock {\em Advances in Mathematics}, 390:107918, 2021.

\bibitem{Handbook_of_Graph_Theory}
Jonathan L. Gross, Jay Yellen, Ping Zhang. 
\newblock {\em Handbook of Graph Theory}. 2nd Edition.
\newblock Chapman and Hall/CRC, 2013.

\bibitem{Kang_Nedela_2021}
M. Kang, R. Nedela, M. {\v{S}}koviera.
\newblock Symmetric maps and their groups.
\newblock {\em Acta Mathematica Slovaca}, 71(2):183--204, 2021.

\bibitem{Krasko_Omelch_4_reg_one_face_maps}
E. Krasko, A. Omelchenko.
\newblock Enumeration of 4-regular one-face maps.
\newblock {\em European Journal of Combinatorics}, 62(5):167--177, 2017.

\bibitem{Krasko_POMI}
E. Krasko.
\newblock Counting unlabelled chord diagrams of maximal genus (in Russian).
\newblock {\em Zapiski Nauchnykh Seminarov POMI}, 464:77--87, 2017.

\bibitem{POMI_Reg_maps_english}
E. Krasko, A. Omelchenko.
\newblock Enumeration of regular maps on surfaces of a given genus.
\newblock {\em Journal of Mathematical Sciences}, 232(1):44--60, 2018.

\bibitem{Torus_Part_I}
E. Krasko, A. Omelchenko.
\newblock Enumeration of $r$-regular maps on the torus. Part I: Rooted maps on
  the torus, the projective plane and the Klein bottle. Sensed maps on the
  torus.
\newblock {\em Discrete Mathematics}, 342(2):584--599, 2019.

\bibitem{Torus_Part_II}
E. Krasko, A. Omelchenko.
\newblock Enumeration of $r$-regular maps on the torus. Part II: Unsensed maps.
\newblock {\em Discrete Mathematics}, 342(2):600--614, 2019.

\bibitem{Unsensed_Maps}
E. Krasko, A. Omelchenko.
\newblock Enumeration of unsensed orientable and non-orientable maps.
\newblock {\em European Journal of Combinatorics}, 86:103093, 2020.

\bibitem{Liskovets_85}
V. Liskovets.
\newblock Enumeration of nonisomorphic planar maps.
\newblock {\em Selecta Math. Sovietica}, 4:303--323, 1985.

\bibitem{Mednykh_Nedela}
A. Mednykh, R. Nedela.
\newblock Enumeration of unrooted maps of a given genus.
\newblock {\em J. Combin. Theory Ser. B}, 96(5):709--729, 2006.

\bibitem{Mednykh_Hypermaps}
A. Mednykh, R. Nedela.
\newblock Enumeration of unrooted hypermaps of a given genus.
\newblock {\em Discrete Mathematics}, 310:518--526, 2010.

\bibitem{Mednykh_Nedela_Stukachev_2023}
A. Mednykh, R. Nedela, A. Stukachev.
\newblock Enumeration of maps with symmetries on surfaces of given genus.
\newblock {\em Mathematical Notes}, 114(6):885--898, 2023.

\bibitem{Tutte_Census}
W. Tutte.
\newblock A census of planar maps.
\newblock {\em Canad. J. Math.}, 15:249--271, 1963.

\bibitem{Tutte_triangulations}
W. Tutte.
\newblock A census of planar triangulations.
\newblock {\em Canad. J. Math.}, 14:21--38, 1962.

\bibitem{Walsh_Lehman}
T. Walsh, A. Lehman.
\newblock Counting rooted maps by genus, I.
\newblock {\em J. Combin. Theory Ser. B}, 13:192--218, 1972.

\bibitem{Walsh_generation}
T. Walsh.
\newblock Space-efficient generation of nonisomorphic maps and hypermaps.
\newblock {\em Journal of Integer Sequences}, 18, 2015.

\end{thebibliography}
\end{document}